# LOCAL LIMIT THEOREMS FOR FINITE AND INFINITE URN MODELS

By Hsien-Kuei Hwang and Svante Janson

*Academia Sinica and Uppsala University*

Local limit theorems are derived for the number of occupied urns in general finite and infinite urn models under the minimum condition that the variance tends to infinity. Our results represent an optimal improvement over previous ones for normal approximation.

**1. Introduction.** A classical theorem of Rényi [29] for the number of empty boxes, denoted by $\mu_0(n,M)$, in a sequence of $n$ random allocations of indistinguishable balls into $M$ boxes with equal probability $1/M$, can be stated as follows: *If the variance of $\mu_0(n,M)$ tends to infinity with $n$, then $\mu_0(n,M)$ is asymptotically normally distributed.* This result, seldom stated in this form in the literature, was proved by Rényi [29] by dissecting the range of $n$ and $M$ into three different ranges, in each of which a different method of proof was employed. Local limit theorems were later studied by Sevast'yanov and Chistyakov [30] in a rather limited range when both ratios of $n/M$ and $M/n$ remain bounded. Kolchin [22] gave a very detailed study on different approximation theorems. For a fairly complete account of this theory, see Kolchin, Sevast'yanov and Chistyakov [23]. Englund [9] later derived an explicit Berry–Esseen bound.

Multinomial extension of the problem was studied by many authors. In this scheme, balls are successively thrown into $M$ boxes, the probability of each ball falling into the $j$th box being $p_j = p_j(M)$, $\sum_{0\le j<M} p_j = 1$. Quine and Robinson [27] showed that *if $p_jM$ is bounded for $j = 0, 1, \ldots, M-1$ and if the variance of $\mu_0(n,M)$ tends to infinity with $n$ ($n$ is the number of allocations), then the distribution of $\mu_0(n,M)$ is asymptotically normal.* They indeed derived a Berry–Esseen bound for the normal approximation of the distribution. Their result remains the strongest of its kind in the literature. Note that the condition $p_jM = O(1)$ for all $j$ is one of the essential conditions needed for proving the asymptotic normality of $\mu_0(n,M)$ in many previous papers on multinomial schemes (see Holst [15], Kolchin, Sevast'yanov and Chistyakov [23]); it implies that the general multinomial scheme studied in the literature is indeed not very far from the equiprobable one. Our major contribution of this paper is to show that this condition can be completely removed. Moreover, under the minimum condition that $\mathrm{Var}(\mu_0(n,M)) \to \infty$, $\mu_0(n,M)$ satisfies a local limit theorem of the form

$$\sup_x \left| \mathbb{P}(\mu_0(n,M) = \mu(n) + x\sigma(n)) - \frac{e^{-x^2/2}}{\sqrt{2\pi}\sigma(n)} \right| = O(\sigma(n)^{-2}), \tag{1.1}$$

for some normalizing constants $\mu(n)$ and $\sigma(n)$, with $\sigma(n) \sim \sqrt{\mathrm{Var}(\mu_0(n,M))}$. Our result is thus, up to the implied constant in the error term, optimum. Moderate and large deviations can also be treated by extending our method of proof, but technicalities will be more involved and the result will be of a less explicit nature; thus we content ourselves with result of the type (1.1).

While the above finite urn schemes have received extensive attention in the literature due to their wide applicability in diverse fields (see Johnson and Kotz [19] and Kotz and Balakrishnan [24]), the model in which $M = \infty$ with $p_j$ fixed was rarely discussed. Bahadur [2], and independently Darling [7], seemed the first to investigate such an urn model. Karlin [20] gave the first systematic study of some basic statistics on this model. His results were then extended by Dutko [8] using the same approach. Dutko showed that in a sequence of $n$ throws, *the distribution of the number of occupied boxes,* $Z_{n,M}$*, is asymptotically normally distributed provided only that its variance tends to infinity with* $n$. (His result was stated in a slightly weaker form.) Note that this result was already stated in the review by Kesten [21] for Darling's [7] paper in *AMS Math Reviews.* For other interesting aspects of $Z_{n,M}$, see the two recent papers [5, 11].

We will derive a local limit theorem for $Z_{n,M}$ of the form (1.1) [with $\mu_0(n,M)$ replaced by $Z_{n,M}$] under the minimum condition that $\mathrm{Var}(Z_{n,M}) \to \infty$, where $M$ is either finite or infinite, and $p_j$ either may depend on $n$ and $M$ or not.

Note that the number of occupied boxes is equivalent to the number of distinct values assumed in a sequence of $n$ i.i.d. (independent and identically distributed) integer-valued random variables. This quantity is an important measure in several problems such as the coupon collector's problem, species-trapping models, birthday paradox, polynomial factorization, statistical linguistics, memory allocation, statistical physics, hashing schemes, and so on. For example, the number of occupied urns under the geometric distribution occurred naturally in at least two different problems in the literature: the depth (the distance of a randomly chosen node to the root) in a class of data structures called Patricia tries (see Rais, Jacquet and Szpankowski [28]) and the number of distinct summands in random integer compositions (see Hitczenko and Louchard [13], Hwang and Yeh [16], Gnedin, Pitman and Yor [12]); see also Prodinger [26] and Janson [18].

Almost all previous approaches rely, explicitly or implicitly, on the widely used Poissonization technique, which is roughly stated as follows. Let $\{a_j\}_j$ be a given sequence such that the Poisson generating function $f(\lambda) := e^{-\lambda}\sum_j a_j\lambda^j/j!$ is an entire function. Then the Poisson heuristic on which the Poissonization procedure relies reads:

(1.2) *if* $f(\lambda)$ *is smooth enough for large* $\lambda$, *then* $a_n \approx f(n)$ $\qquad (n\to\infty)$.

Such a heuristic, guided by the underlying normal approximation to the Poisson distribution, can usually be justified by suitable real or complex analysis. As is often the case, it is the verification of the smoothness (or regularity) property of $f(\lambda)$ that is the hard part of the heuristic and for which technical conditions are usually introduced. The heuristic appeared in different guises in diverse contexts such as Borel summability and Tauberian theorems; it can at least be traced back to Ramanujan's Notebooks; see the book by Berndt [4], pages 57–66, for more details, Aldous [1] and the survey paper by Jacquet and Szpankowski [17] for thorough discussions.

To obtain our local limit theorems, we apply instead the two-dimensional saddle-point method, which is in essence the most straightforward one and may be regarded as an extension of the Poisson heuristic; see also Remark 3.2. The approach we use can be extended in a few lines: moderate and large deviations of $Z_{n,M}$, consideration of other statistics such as urns with a given number of balls, weighted coverage, goodness-of-fit tests, etc.

This paper is organized as follows. We first state our main results on local limit theorems in the next section. In Section 3 the case of a Poisson number of balls is considered, and we introduce the Poisson generating function that is central to our proofs. Asymptotics of mean and variance are derived in Section 4. Sections 5–7 give the proofs of the main results. Discrete limit laws are derived in Section 8. We conclude this paper with some properties of infinite urn models.

*Notation.* The generic symbols $C, C_1, C_2, \dots$ and $c_1, c_2, \dots$ will always denote some positive absolute constants; they can be replaced by explicit numerical values if desired, but we avoid this for simplicity of presentation. Similarly, the implicit constants in the $O$- and $\asymp$-symbols are absolute constants, where the symbol $A \asymp B$ means that $c \le A/B \le C$ for some constants $c$ and $C$.

**2. Results.** Let $X_1, X_2, \dots, X_n$ be a sequence of i.i.d. random variables with a discrete distribution $F$. Let $Z = Z_{n,F}$ denote the number of distinct values assumed by $X_1, \dots, X_n$.

Let $\mathcal{J}$ be the set (finite or infinite) of possible values of $X_i$, and let the distribution $F$ be given by

$$\mathbb{P}(X_i = j) = p_j \qquad (j \in \mathcal{J}),$$

where $\sum_j p_j = 1$. *Here and throughout this paper*, sums of the form $\sum_j$ are taken to be $\sum_{j\in\mathcal{J}}$ unless otherwise specified; similarly $\prod_j = \prod_{j\in\mathcal{J}}$.

Alternatively, $Z$ counts the number of occupied urns in an urn scheme where $n$ balls are thrown independently and each ball has the same probability $p_j$ of falling into urn $j$, $j \in \mathcal{J}$. Note that we allow $p_j = 0$ for some $j$, although such elements may be freely added to or deleted from $\mathcal{J}$ without changing $Z$.

We now state our results. Proofs are given in Sections 4–7.

THEOREM 2.1. *If* $\mathrm{Var}(Z_{n,F}) \neq 0$*, then the local limit theorem*

$$\sup_{-\infty<x<\infty}\left|\mathbb{P}(Z_{n,F} = m) - \frac{e^{-x^2/2}}{\sqrt{2\pi\,\mathrm{Var}(Z_{n,F})}}\right| \le \frac{C}{\mathrm{Var}(Z_{n,F})} \tag{2.1}$$

*holds uniformly for all* $n \in \mathbb{N}$ *and* $F$*, where* $m = \lfloor \mathbb{E}(Z_{n,F}) + x\sqrt{\mathrm{Var}(Z_{n,F})}\rfloor$.

The trivial case $\mathrm{Var}(Z_{n,F}) = 0$ occurs if and only if $n = 0$, $n = 1$, or $F$ is a one-point distribution; in these cases $Z = 0$, 1 and 1, respectively.

REMARK 2.1 (*Discrete distributions vs. continuous distributions*). The assumption that the distribution $F$ is discrete is not necessary. If $F$ is continuous, then $Z = n$ a.s., another trivial case with $\mathrm{Var}(Z) = 0$. If $F$ has both a discrete and a continuous part, then Theorem 2.1 still holds. To see this, assume that $F$ has a continuous part with total mass $\rho$, and let $F_M$ be a discrete distribution that has the same atoms as $F$ together with $M$ new atoms $j$, each with $p_j = \rho/M$. We can now apply Theorem 2.1 to $F_M$, and it is easily seen that if we let $M \to \infty$ (with $n$ fixed), then all quantities in (2.1) for $F_M$ converge to the corresponding quantities for $F$; thus (2.1) holds for $F$ also.

Similarly, the result below holds for general distributions with minor modifications in the formulas (2.2)–(2.5) for mean and variance. We omit the details.

REMARK 2.2 (*Finite urns vs. infinite urns*). By a suitable truncation, it suffices to prove the results for a finite set $\mathcal{J}$. This has the technical advantage that we do not need to address the convergence of the sums and products involved, which is, however, relatively easily checked. Indeed, without loss of generality, we may assume that $\mathcal{J}$ is the set of nonnegative integers; then we replace $X_i$ by $X_i \wedge M$, and let $F_M$ be the corresponding distribution (i.e., $F$ truncated at $M$). It follows that if the result holds for each $F_M$, it also holds for $F$, by letting $M \to \infty$.

Exact formulas for $\mathbb{E}(Z_{n,F})$ and $\operatorname{Var}(Z_{n,F})$ are given in (4.1) and (4.2) below. However, these formulas are rather complicated; thus we first derive simpler approximations to these quantities.

We define, for $x \geq 0$ (and, more generally, for any complex $x$ with $\Re x \geq 0$),

$$\mu_F(x) := \sum_j (1 - e^{-p_j x}), \tag{2.2}$$

$$v_F(x) := \sum_j e^{-p_j x}(1 - e^{-p_j x}), \tag{2.3}$$

$$u_F(x) := \sum_j p_j x e^{-p_j x}, \tag{2.4}$$

$$\sigma_F^2(x) := v_F(x) - \frac{u_F(x)^2}{x} \tag{2.5}$$

(with $\sigma_F^2(0) := 0$) and, for later use,

$$\tilde{v}_F(x) := x + v_F(x) - 2u_F(x). \tag{2.6}$$

We will see in Section 3 that $\mu_F(x)$ and $v_F(x)$ are the mean and variance of $Z$ if the fixed number $n$ of variables (balls) is replaced by a Poisson number with mean $x \geq 0$, and that $u_F(x)$, $\sigma_F^2(x)$ and $\tilde{v}_F(x)$ also have simple interpretations in terms of this Poissonized version. Noting that

$$\tfrac{1}{2} \sum_{i,j \in \mathcal{J}} p_i p_j (e^{-p_i x} - e^{-p_j x})^2 = \sum_i p_i e^{-2p_i x} - \Big( \sum_i p_i e^{-p_i x} \Big)^2,$$

we obtain the following alternative formula.

Proposition 2.2.

$$\sigma_F^2(x) = \sum_j e^{-p_j x}(1 - (1 + p_j x) e^{-p_j x}) + \frac{x}{2} \sum_{i,j} p_i p_j (e^{-p_i x} - e^{-p_j x})^2. \tag{2.7}$$

All terms in the sums in (2.7) are nonnegative for $x > 0$. Hence, $\sigma_F^2(x) > 0$ for any $F$ and all $x > 0$.

Theorem 2.3. *The mean and the variance of $Z_{n,F}$ satisfy*

$$\mathbb{E}(Z_{n,F}) = \mu_F(n) + O(1), \tag{2.8}$$

$$\operatorname{Var}(Z_{n,F}) = \sigma_F^2(n) + O(1). \tag{2.9}$$

The $O(1)$-terms in (2.8) and (2.9) are in some cases $o(1)$, as we will see later.

We can thus replace the exact mean and the exact variance in Theorem 2.1 by their asymptotic approximations.

THEOREM 2.4. *If $\sigma_F(n) \neq 0$, then uniformly for all $n \geq 1$ and $F$*

$$\sup_{-\infty<x<\infty}\left|\mathbb{P}(Z_{n,F} = \lfloor \mu_F(n) + x\sigma_F(n)\rfloor) - \frac{e^{-x^2/2}}{\sqrt{2\pi}\sigma_F(n)}\right| \leq \frac{C_1}{\sigma_F^2(n)}.$$

These results are stated as approximation results. If we consider a sequence of such variables $Z_{n,F}$, by letting $n \to \infty$ and varying $F$, assuming only that $\operatorname{Var}(Z_{n,F}) \to \infty$, Theorems 2.1 and 2.4 can be interpreted as local central limit theorems. The corresponding central limit theorem, with (the generally weaker) convergence in distribution, can be stated as follows.

COROLLARY 2.5. *Consider a sequence $(n_\nu, F_\nu)_\nu$ of integers $n_\nu$ and distributions $F_\nu$. Then the following statements are equivalent, with $Z_\nu := Z_{n_\nu,F_\nu}$ and $\sigma_\nu^2 := \operatorname{Var}(Z_\nu)$:*

(i) $\sigma_\nu^2 \to \infty$;
(ii) $\sigma_{F_\nu}^2(n_\nu) \to \infty$;
(iii) $(Z_\nu - \mathbb{E}(Z_\nu))/\sigma_\nu \xrightarrow{\mathrm{d}} N(0,1)$;
(iv) $(Z_\nu - \mu_{F_\nu}(n_\nu))/\sigma_{F_\nu}(n_\nu) \xrightarrow{\mathrm{d}} N(0,1)$;
(v) $(Z_\nu - \alpha_\nu)/\beta_\nu \xrightarrow{\mathrm{d}} N(0,1)$ *for some sequences $\alpha_\nu$ and $\beta_\nu$ with $\beta_\nu > 0$.*

These theorems cover many results in the literature as special cases.

From now on, the distribution $F$ will be fixed, and we will generally drop the subscript $F$ from the notation.

For our method of proof, we consider two (technical) cases:

(i) $\sum_{p_j n \leq 1} p_j \geq 1/2$, meaning roughly that asymptotics of $Z_n$ is dominated by small $p_j$.

(ii) $\sum_{p_j n > 1} p_j \geq 1/2$, meaning roughly that asymptotics of $Z_n$ is dominated by large $p_j$.

Obviously, at least one of these cases will hold. Here, the value $1/2$ is not essential and can be changed to any small positive constant (with consequent changes in the values of some of the unspecified constants below); similarly, the cut-off at $p_j n = 1$ is chosen for technical convenience. We will work with $Z_n$ in case (ii) and with $\tilde{Z}_n := n - Z_n$ in case (i); it will turn out that Poissonization then works well in both cases.

REMARK 2.3 (*Exact distribution*). It is easy to find the exact distribution of $Z_{n,F}$. Indeed, assuming as we may that the set $\mathcal{J}$ is ordered,

$$\mathbb{P}(Z_{n,F}=m)=\sum_{\substack{h_1+\cdots+h_m=n\\ h_j\geq 1}}\sum_{j_1<\cdots<j_m}\binom{n}{h_1,\dots,h_m}p_{j_1}^{h_1}\cdots p_{j_m}^{h_m}. \tag{2.10}$$

This expression explains why such urn schemes are called multinomial allocations. However, it will not be used in this paper.

**3. Poissonization.** We consider first the mean and the variance of the number of occupied urns or the number of distinct values when the number of balls or variables have a Poisson distribution.

Recall that $Z_n$ is the number of occupied urns when we throw $n$ balls. Then $\tilde{Z}_n := n - Z_n$ represents the number of balls that land in a nonempty urn.

Consider now instead the case when the number $N$ of balls is Poisson distributed. Let $Z(\lambda)$ denote the number of occupied balls with $N = N(\lambda) \sim \mathrm{Po}(\lambda)$ balls; let $\tilde{Z}(\lambda) := N(\lambda) - Z(\lambda)$ be the number of balls that land in an occupied urn.

REMARK 3.1 (*A coupling*). We may define $Z_n$ and $Z(\lambda)$ for all $n \geq 0$ and $\lambda \geq 0$ simultaneously (for a given $F$) by throwing balls at times given by a Poisson process with intensity 1. We let $Z(\lambda)$ be the number of occupied urns at time $\lambda \geq 0$, when $\mathrm{Po}(\lambda)$ balls have been thrown, and let $Z_n$ be the number of occupied urns when the $n$th ball has been thrown. This defines the various variables simultaneously, with both $Z_n$ and $\tilde{Z}_n := n - Z_n$ increasing in $n$ and both $Z(\lambda)$ and $\tilde{Z}(\lambda) := N(\lambda) - Z(\lambda)$ increasing in $\lambda$, where $N(\lambda) \sim \mathrm{Po}(\lambda)$ is the number of balls thrown at time $\lambda$.

Let $U_j$ be the number of balls in urn $j$. Then $N = \sum_j U_j$,

$$Z = \sum_j \mathbf{1}_{\{U_j \geq 1\}}, \tag{3.1}$$

where $\mathbf{1}_{\mathcal{A}}$ denotes the indicator function of the event $\mathcal{A}$, and

$$\tilde{Z} = \sum_j U_j - Z = \sum_j (U_j - \mathbf{1}_{\{U_j \geq 1\}}) = \sum_j \tilde{U}_j,$$

where $\tilde{U}_j := U_j - \mathbf{1}_{\{U_j \geq 1\}}$.

In the Poisson case, the random variables $U_j$ are independent and Poisson distributed, with $U_j \sim \mathrm{Po}(p_j\lambda)$. By (3.1), $Z(\lambda) = \sum_j I_j$, where $I_j :=$ $\mathbf{1}_{\{U_j\geq 1\}} \sim \mathrm{Be}(1-e^{-p_j\lambda})$ are independent Bernoulli random variables. It follows that

$$\mathbb{E}(Z(\lambda)) = \sum_j \mathbb{E}(I_j) = \sum_j (1-e^{-p_j\lambda}) = \mu(\lambda),$$
$$\mathrm{Var}(Z(\lambda)) = \sum_j \mathrm{Var}(I_j) = \sum_j e^{-p_j\lambda}(1-e^{-p_j\lambda}) = v(\lambda).$$

Similarly, $\tilde{Z}(\lambda) = \sum_j \tilde{U}_j$, where $\tilde{U}_j = (U_j + \mathbf{1}_{\{U_j=0\}} - 1)$ are independent. We have

$$\mathbb{E}(\tilde{U}_j) = p_j\lambda + e^{-p_j\lambda} - 1,$$
$$\begin{aligned}\mathrm{Var}(\tilde{U}_j) &= \mathrm{Var}(U_j) + \mathrm{Var}(\mathbf{1}_{\{U_j=0\}}) + 2\,\mathrm{Cov}(U_j, \mathbf{1}_{\{U_j=0\}}) \\ &= p_j\lambda + e^{-p_j\lambda}(1-e^{-p_j\lambda}) - 2p_j\lambda e^{-p_j\lambda}.\end{aligned}$$

Accordingly [see (2.2)–(2.6)],

$$\mathbb{E}(\tilde{Z}(\lambda)) = \sum_j (p_j\lambda - (1-e^{-p_j\lambda})) = \lambda - \mu(\lambda),$$
$$\begin{aligned}\mathrm{Var}(\tilde{Z}(\lambda)) &= \sum_j (p_j\lambda + e^{-p_j\lambda}(1-e^{-p_j\lambda}) - 2p_j\lambda e^{-p_j\lambda}) \\ &= \lambda + v(\lambda) - 2u(\lambda) = \tilde{v}(\lambda).\end{aligned} \tag{3.2}$$

Remark 3.2 (*A connection between the two cases* $\sum_{p_jn\leq 1} p_j \geq \frac{1}{2}$ *and* $\sum_{p_jn>1} p_j \geq \frac{1}{2}$). We have seen that $\mu(\lambda)$ and $v(\lambda)$ are the mean and variance of $Z$ in the Poisson case. Similarly, it is easy to see that $\mathrm{Cov}(Z,N) = u(\lambda)$, and thus

$$\sigma^2(\lambda) = \mathrm{Var}(Z) - \frac{\mathrm{Cov}(Z,N)^2}{\mathrm{Var}\,N}, \tag{3.3}$$

which can be interpreted to be the smallest variance of a linear combination $Z - \alpha N$ with $\alpha \in \mathbb{R}$, that is, $\sigma^2(\lambda) = \mathrm{Var}(Z - \alpha_0 N)$, where $Z - \alpha_0 N$ is optimal in this sense, which, on the other hand, is also determined by $\mathrm{Cov}(Z - \alpha_0 N, N) = 0$. (These are standard equations in linear regression, where (3.3) is the residual variance, and $\alpha_0 = \mathrm{Cov}(Z,N)/\mathrm{Var}(N) \in [0,1]$ because $\mathrm{Cov}(Z,N) = u(\lambda) \leq \lambda = \mathrm{Var}(N)$.)

Our method of proof is based on analyzing the Poisson generating function $P$ defined below. This can be regarded as an analytical Poissonization, and it is, at least heuristically, strongly related to replacing $Z_n$ by the Poissonized $Z(\lambda)$ and then compensating for the randomness in $N$, the number of balls, in order to derive results for $Z_n$. It is then natural to consider the projection $Z-\alpha_0 N$, which eliminates the first-order (linear) fluctuations in $Z$ due to the randomness in $N$. Theorem 2.3 says that, with $\lambda=n$, this projection has almost the same variance as $Z_n$, which indicates that this projection (plus the constant $\alpha_0 n$) is a good approximation to $Z_n$. Moreover, as we will see in Proposition 4.3 below, the smallest variance $\sigma^2(\lambda)$ of a linear combination $Z-\alpha N$ is attained within a constant factor by one of the choices $\alpha=0$ and $\alpha=1$, which gives $Z$ and $Z-N=-\tilde{Z}$, respectively. Indeed, the arguments below can be interpreted as considering these two choices only. (It is likely that one could use similar arguments corresponding to the optimal projection $Z-\alpha_0 N$, without splitting our analysis below into two cases.)

We define, for complex $z$ and $y$, $P(z,y)$ to be the exponential generating function of $\mathbb{E}(y^{Z_n})$ given by

$$P(z,y) := \sum_{n\geq 0} \frac{z^n}{n!}\mathbb{E}(y^{Z_n}) = \sum_{n,m\geq 0} \frac{z^n y^m}{n!}\mathbb{P}(Z_n=m); \tag{3.4}$$

we further define the Poisson generating function $Q(z,y) := e^{-z}P(z,y)$. Note that for $\lambda\geq 0$, $Q(\lambda,y)=\mathbb{E}(y^{Z(\lambda)})$ is the probability generating function of $Z(\lambda)$. It follows immediately from (3.1) for $z\geq 0$, and for general complex $z$ by analytic continuation, that

$$Q(z,y)=\prod_j (1+(y-1)(1-e^{-p_j z})),$$

and thus, using $\sum_j p_j=1$,

$$P(z,y)=e^z Q(z,y)=\prod_j (1+y(e^{p_j z}-1)). \tag{3.5}$$

This also follows easily from (2.10); see also Karlin [20], Johnson and Kotz [19], Kolchin, Sevast'yanov and Chistyakov [23], Flajolet, Gardy and Thimonier [10] for different derivations.

Note also that the probability generating function of $\tilde{Z}(\lambda)=N-Z(\lambda)$ is given by

$$\mathbb{E}(y^{\tilde{Z}(\lambda)}) = \sum_{n\geq 0}\frac{\lambda^n}{n!}e^{-\lambda}\mathbb{E}(y^{n-Z_n}) = e^{-\lambda}P(y\lambda, y^{-1}). \tag{3.6}$$

According to the Poisson heuristic (1.2), if $Q$ were smooth enough, then we would have

$$\mathbb{E}(y^{Z_n})\approx Q(n,y) \qquad (n\to\infty),$$

and the asymptotic normality of $Z_n$ would then follow from Taylor expansion of the cumulant generating function

$$\log \mathbb{E}(e^{Z_n s}) \approx s\sum_j (1-e^{-p_j n}) + \frac{s^2}{2}\sum_j e^{-p_j n}(1-e^{-p_j n})+\cdots,$$

provided that the second sum tends to infinity and the error term becomes small after normalization. However, the general situation here turns out to be more complicated. First, the variance of $Z_n$ is not necessarily of the same order as the second sum. Second, the growth order of $Q(z,y)$ is not necessarily polynomial in $|z|$; for example, $Q(z,0)=e^{-z}$. Thus more refined arguments are required to properly justify the (implicit) underlying Poisson heuristic (1.2).

**4. Mean and variance of $Z$.** We prove in this section the estimates (2.8) and (2.9) for the mean and variance of $Z_n$, and some related estimates.

Proof of Theorem 2.3. By straightforward calculations, (3.1) leads to

$$\mathbb{E}(Z_n)=\sum_j (1-(1-p_j)^n), \tag{4.1}$$

$$\begin{aligned}\mathrm{Var}(Z_n)&=\sum_j (1-p_j)^n(1-(1-p_j)^n)\\ &\quad+\sum_{i\neq j}((1-p_i-p_j)^n-(1-p_i)^n(1-p_j)^n).\end{aligned} \tag{4.2}$$

Now for $p\in[0,1]$, we have

$$\begin{aligned}0&\le e^{-pn}-(1-p)^n\le ne^{-p(n-1)}(e^{-p}-1+p)\\ &=O(p^2ne^{-pn})=O(p),\end{aligned} \tag{4.3}$$

and thus (2.8) follows from (4.1).

Similarly, by $(1-x)^n=1-nx+O(n^2x^2)$ for $0\le x\le 1$, we have, for $n\ge 2$,

$$\begin{aligned}&(1-p_i)^n(1-p_j)^n-(1-p_i-p_j)^n\\ &\quad=(1-p_i)^n(1-p_j)^n\left(1-\left(1-\frac{p_ip_j}{(1-p_i)(1-p_j)}\right)^n\right)\\ &\quad=(1-p_i)^n(1-p_j)^n\left(\frac{p_ip_jn}{(1-p_i)(1-p_j)}+O\left(\frac{p_i^2p_j^2n^2}{(1-p_i)^2(1-p_j)^2}\right)\right)\\ &\quad=p_ip_jn(1-p_i)^{n-1}(1-p_j)^{n-1}+O((p_ip_jn)^2(1-p_i)^{n-2}(1-p_j)^{n-2})\\ &\quad=p_ip_jne^{-p_in-p_jn}+O(p_ip_j),\end{aligned}$$

since, by (4.3),

$$\begin{aligned}e^{-pn}-(1-p)^{n-1}&=e^{-pn}-(1-p)^n+p(1-p)^{n-1}\\ &=O(p^2ne^{-pn}+pe^{-pn})=O(1/n),\end{aligned}$$

and

$$(p_i p_j n)^2 (1-p_i)^{n-2}(1-p_j)^{n-2} = O(p_i p_j \cdot p_i n e^{-p_i n} \cdot p_j n e^{-p_j n}) = O(p_i p_j).$$

Hence (4.2) yields

$$\begin{aligned}
\operatorname{Var}(Z_n) &= \sum_j (e^{-p_j n} + O(p_j) - e^{-2p_j n} + O(p_j)) \\
&\quad - \sum_{i\neq j} (p_i p_j n e^{-p_i n - p_j n} + O(p_i p_j)) \\
&= \sum_j (e^{-p_j n} - e^{-2p_j n}) - \sum_{i,j} p_i p_j n e^{-p_i n - p_j n} + \sum_i p_i^2 n e^{-2p_i n} + O(1) \\
&= \sum_j (e^{-p_j n} - e^{-2p_j n}) - n\left(\sum_i p_i e^{-p_i n}\right)^2 + O(1),
\end{aligned}$$

which proves (2.9). □

We proceed to some estimates of $v(x)$, $\tilde{v}(x)$ and $\sigma^2(x)$, which roughly indicate why we need to separate into the cases $p_j n > 1$ and $p_j n \le 1$ in our manipulations of sums.

LEMMA 4.1. *For $x \ge 0$, we have*

$$v(x) = \operatorname{Var}(Z(x)) \asymp \sum_{p_j x \le 1} p_j x + \sum_{p_j x > 1} e^{-p_j x}, \tag{4.4}$$

$$\tilde{v}(x) = \operatorname{Var}(\tilde{Z}(x)) \asymp \sum_{p_j x \le 1} (p_j x)^2 + \sum_{p_j x > 1} p_j x. \tag{4.5}$$

PROOF. This follows from the definitions (2.3) and (2.6) [see also (3.2)] and the asymptotics

$$e^{-x}(1-e^{-x}) \sim \begin{cases} x, & \text{if } x \to 0, \\ e^{-x}, & \text{if } x \to \infty, \end{cases}$$

and

$$x + e^{-x}(1 - e^{-x} - 2x) \sim \begin{cases} x^2/2, & \text{if } x \to 0, \\ x, & \text{if } x \to \infty. \end{cases}$$

□

LEMMA 4.2. *For all $x \ge 0$*

$$\sigma^2(x) \le \sum_{p_j x \le 1} (p_j x)^2 + \sum_{p_j x > 1} p_j x \tag{4.6}$$

*and*

$$\sigma^2(x) \ge c_1 \left( \sum_{p_j x \le 1} (p_j x)^2 + \sum_{p_j x > 1} e^{-p_j x} \right). \tag{4.7}$$

PROOF. Since $v(x) \le u(x) \le x$, we have

$$\sigma^2(x) = v(x) - \frac{u^2(x)}{x} \le \frac{u(x)}{x}(x - u(x)) \le x - u(x) = \sum_j p_j x(1 - e^{-p_j x}),$$

from which the upper bound in (4.6) follows.

On the other hand, by Proposition 2.2, $\sigma^2(x) \ge \sum_j e^{-p_j x}(1-(1+p_j x)e^{-p_j x})$, which yields (4.7) by the elementary inequality $1-(1+x)e^{-x} \ge c_2 \min\{1, x^2\}$. □

Note that by the inequality $1-(1+x)e^{-x} \ge x^2 e^{-x}/2$, we also have $\sigma^2(x) \ge \frac{1}{2}\sum_j (p_j x)^2 e^{-2p_j x}$.

The following result, based on the estimates we just derived, is crucial for the development of our arguments.

PROPOSITION 4.3. $\sigma^2(x) \asymp \min(v(x), \tilde{v}(x))$. *More precisely:*

(i) *if* $\sum_{p_j x \le 1} p_j \ge 1/2$, *then*

$$\sigma^2(x) \asymp \tilde{v}(x);$$

(ii) *if* $\sum_{p_j x > 1} p_j \ge 1/2$, *then*

$$\sigma^2(x) \asymp v(x).$$

PROOF. The upper bounds are immediate: $\sigma^2(x) \le v(x)$ by definition, while $\sigma^2(x) = O(\tilde{v}(x))$ by (4.5) and (4.6). Alternatively, as pointed out by one of the referees, the upper bounds also follow from Remark 3.2 and $\operatorname{Var}(Z) = v(\lambda)$, $\operatorname{Var}(\tilde{Z}) = \operatorname{Var}(Z - N) = \tilde{v}(\lambda)$.

For the lower bounds, we treat the two cases separately.

Case (i): $\sum_{p_j x \le 1} p_j \ge 1/2$. We have

$$\begin{aligned} x\sum_{i,j} p_i p_j (e^{-p_i x} - e^{-p_j x})^2 &\ge x \sum_{p_i x \le 1} \sum_{p_j x \ge 2} p_i p_j (e^{-1} - e^{-2})^2 \\ &\ge \frac{(e^{-1} - e^{-2})^2}{2} \sum_{p_j x \ge 2} p_j x. \end{aligned}$$

Thus, using Proposition 2.2,

$$\sum_{p_j x \ge 2} p_j x = O(\sigma^2(x)).$$

Moreover, by (4.7),

$$\sum_{p_j x \le 1} (p_j x)^2 + \sum_{1 < p_j x \le 2} p_j x = O(\sigma^2(x)).$$

Hence, by (4.5), $\tilde{v}(x) = O(\sigma^2(x))$.

Case (ii): $\sum_{p_j x > 1} p_j \geq 1/2$. First, we have by Proposition 2.2

$$\begin{aligned}\sum_{p_i x \leq 1/2} p_i x &\leq 2x \sum_{p_i x \leq 1/2} \sum_{p_j x > 1} p_i p_j \\ &= O\Bigg( x \sum_{i,j} p_i p_j (e^{-p_i x} - e^{-p_j x})^2 \Bigg) \\ &= O(\sigma^2(x)).\end{aligned}$$

Furthermore, by (4.7),

$$\sum_{1/2 < p_j x \leq 1} p_j x + \sum_{p_j x > 1} e^{-p_j x} = O(\sigma^2(x)). \tag{4.8}$$

Thus, by (4.4), $v(x) = O(\sigma^2(x))$. □

Remark 4.1 (*An interesting estimate*). It follows from Lemma 4.1 that in case (i), $v(x) \asymp x$, and in case (ii), $\tilde{v}(x) \asymp x$. Hence, $\max(v(x), \tilde{v}(x)) \asymp x$.

## 5. Local limit theorem when $\sum_{p_j n > 1} p_j \geq 1/2$.

We prove Theorem 2.4 in this section when $\sum_{p_j n > 1} p_j \geq 1/2$. Our starting point is the integral representation

$$\mathbb{P}(Z_n = m) = \frac{n! n^{-n}}{(2\pi)^2} \int_{-\pi}^{\pi} \int_{-\pi}^{\pi} e^{-\mathrm{i}m\varphi - \mathrm{i}n\theta} P(ne^{\mathrm{i}\theta}, e^{\mathrm{i}\varphi}) \, d\theta \, d\varphi, \tag{5.1}$$

which follows from (3.4) by standard coefficient extraction.

Our strategy is to apply the two-dimensional saddle-point method. More precisely, we split the integration ranges of the double integral into three parts:

$$\iint_{\substack{|\theta| \leq \theta_0 \\ |\varphi| \leq \varphi_0}} + \iint_{\substack{|\theta| \leq \theta_0 \\ \varphi_0 < |\varphi| \leq \pi}} + \iint_{\substack{\theta_0 < |\theta| \leq \pi \\ |\varphi| \leq \pi}}, \tag{5.2}$$

where $\theta_0$ and $\varphi_0$ are usually so chosen that they satisfy the conditions for the saddle-point method:

$$n\theta_0^2 \to \infty, \qquad n\theta_0^3 \to 0 \qquad \text{and} \qquad \sigma(n)^2 \varphi_0^2 \to \infty, \qquad \sigma(n)^2 \varphi_0^3 \to 0.$$

For technical convenience, we will instead choose $\theta_0 := n^{-1/2} \sigma(n)^{1/3} \leq n^{-1/3}$ and $\varphi_0 := \sigma(n)^{-2/3}$, and the usual saddle-point method will require only minor modifications.

We show that the main contribution to $\mathbb{P}(Z_n = m)$ comes from the first double integral in (5.2), the other two being asymptotically negligible. As is often the case, the hard part of the proof is to prove the smallness of $e^{-n} |P(ne^{\mathrm{i}\theta}, e^{\mathrm{i}\varphi})|$ when at least one of $\{\theta, \varphi\}$ is away from zero. Note that $P(n, 1) = e^n$.

5.1. *Estimates for* $|P(z, e^{i\varphi})|$. We derive in this subsection two major estimates for $|P(ne^{i\theta}, e^{i\varphi})|$ under the assumption $\sum_{p_j n>1} p_j \ge 1/2$. The corresponding estimates for the case $\sum_{p_j n\le 1} p_j \ge 1/2$ will be given in the next section.

LEMMA 5.1. *Let* $z = re^{i\theta}$, $r \ge 0$ *and* $\theta \in \mathbb{R}$. *Then:*

$$|e^z - 1| \le (e^r - 1)e^{-r(1-\cos\theta)/2}.$$

PROOF. We have

$$|e^z - 1| = 2|e^{z/2}||\sinh(z/2)| \le 2e^{(r/2)\cos\theta}\sinh(r/2) = e^{(r/2)\cos\theta}(e^{r/2} - e^{-r/2}),$$

from which the result follows. □

LEMMA 5.2. *If* $r \ge 1$ *and* $|\theta| \le \pi$, *then*

$$1 + |e^{re^{i\theta}} - 1| \le e^{r - c_3 r\theta^2}. \tag{5.3}$$

PROOF. By Lemma 5.1

$$1 + |e^{re^{i\theta}} - 1| \le 1 + (e^r - 1)e^{-r(1-\cos\theta)/2} \le 1 + (e^r - 1)e^{-c_4 r\theta^2}, \tag{5.4}$$

for $|\theta| \le \pi$, where we can take $c_4 = 1/\pi^2$ by the inequality $1 - \cos\theta \ge 2\theta^2/\pi^2$ for $|\theta| \le \pi$. Define $c_3 := c_4/2$. By the inequalities

$$e^{c_3 r\theta^2} + 1 \le e^{r/2} + 1 \le e^r,$$

we have

$$1 - e^{-2c_3 r\theta^2} = (e^{c_3 r\theta^2} + 1)e^{-c_3 r\theta^2}(1 - e^{-c_3 r\theta^2}) \le e^r(e^{-c_3 r\theta^2} - e^{-2c_3 r\theta^2}).$$

The result (5.3) follows from this and (5.4). □

PROPOSITION 5.3. *Assume that* $\sum_{p_j n>1} p_j \ge 1/2$. *Then the inequality*

$$|P(ne^{i\theta}, e^{i\varphi})| \le e^{n - c_5 n\theta^2}$$

*holds uniformly for* $|\theta| \le \pi$ *and* $-\infty < \varphi < \infty$, *where* $c_5 = c_3/2$.

PROOF. By (3.5), Lemma 5.2 and the simple estimate $1 + |e^{re^{i\theta}} - 1| \le e^r$ (e.g., by Lemma 5.1),

$$\begin{aligned}|P(ne^{i\theta}, e^{i\varphi})| &\le \prod_j (1 + |e^{p_j n e^{i\theta}} - 1|)\\ &\le \Bigg(\prod_{p_j n\le 1} e^{p_j n}\Bigg)\Bigg(\prod_{p_j n>1} e^{p_j n - c_3 p_j n\theta^2}\Bigg)\\ &= \exp\Bigg(n - c_3\theta^2 \sum_{p_j n>1} p_j n\Bigg).\end{aligned}$$

□

Since a more detailed estimate for $|P(ne^{\mathrm{i}\theta}, e^{\mathrm{i}\varphi})|$ and a local expansion of $P(ne^{\mathrm{i}\theta}, e^{\mathrm{i}\varphi})$ (for small $\theta$ and $\varphi$) involve several sums related to $u(x), v(x), \tilde{v}(x)$ and $\sigma^2(x)$, we now derive a few simple estimates for and relationships between them.

LEMMA 5.4. *For $x \ge 0$,*

$$u(x) = O(x \wedge (v(x) + 1)(1 + \log_+ x)).$$

PROOF. The upper bound $x$ follows easily by the inequality $e^{-x} \le 1$. For the other upper bound, we may assume $x \ge 2$. Then, by (4.4),

$$\sum_{p_j x \le \log x} p_j x e^{-p_j x} \le \sum_{p_j x \le 1} p_j x + \log x \sum_{p_j x > 1} e^{-p_j x} = O(v(x) \log x), \tag{5.5}$$

while trivially $\sum_{p_j x > \log x} p_j x e^{-p_j x} \le \sum_j p_j = 1$. □

LEMMA 5.5. *Let $x \ge 1$ and $0 \le \delta \le 1/2$. Then:*

(i) *$v(x) = O(v((1-\delta)x))$ and $v((1-\delta)x) = O(x^{2\delta} v(x) + 1)$;*
(ii) *$\tilde{v}(x) \asymp \tilde{v}((1-\delta)x)$;*
(iii) *$\sigma^2(x) = O(\sigma^2((1-\delta)x))$;*
(iv) *if, furthermore, $\delta \le x^{-1/3}$, then $v((1-\delta)x) = O(v(x)+1)$ and $\sigma^2((1-\delta)x) = O(\sigma^2(x) + 1)$.*

PROOF. (i). We use (4.4) for both $x$ and $(1-\delta)x$; note that we can split the sum according to $p_j x \le 1$ and $p_j x > 1$ for $(1-\delta)x$, too. The first estimate then is obvious. For the second we find, assuming as we may $x \ge 2$,

$$\sum_{p_j x \le 1} p_j (1-\delta) x \le \sum_{p_j x \le 1} p_j x = O(v(x)),$$

$$\sum_{1 < p_j x \le 2 \log x} e^{-p_j (1-\delta) x} \le e^{2\delta \log x} \sum_{1 < p_j x \le 2 \log x} e^{-p_j x} \le C_2 x^{2\delta} v(x),$$

$$\sum_{p_j x > 2 \log x} e^{-p_j (1-\delta) x} \le \sum_{p_j x > 2 \log x} e^{-\log x} \le 1,$$

since there are at most $x/\log x$ terms in the last sum.

(ii). Immediate from (4.5).

(iii) and (iv). Follow from (i) and (ii) together with Proposition 4.3. □

We now refine Proposition 5.3 and obtain a decrease of $|P(ne^{\mathrm{i}\theta}, e^{\mathrm{i}\varphi})|$ in both $\theta$ and $\varphi$. (We are grateful to one of the referees for improving our previous version.)

Proposition 5.6. *Assume that* $\sum_{p_jn>1} p_j \ge 1/2$. *Then uniformly for* $|\theta| \le \pi$ *and* $|\varphi| \le \pi$, *provided* $v(n) \ge 1$,

$$|P(ne^{\mathrm{i}\theta}, e^{\mathrm{i}\varphi})| \le \exp(n - c_6 n\theta^2 - c_7\sigma^2(n)\varphi^2). \tag{5.6}$$

Proof. Let $z := ne^{\mathrm{i}\theta} = \xi + \mathrm{i}\eta$, where $\xi := n\cos\theta$ and $\eta := n\sin\theta$. Assume first that $|\theta| \le \pi/4$; then $|\eta| \le \xi$ and $\xi \ge n/2$. Thus $n = O(\xi)$ and, by Lemma 5.5(i), $1 \le v(n) = O(v(\xi))$. By explicit calculation

$$\begin{aligned}
&|1 + e^{\mathrm{i}\varphi}(e^{\xi+\mathrm{i}\eta} - 1)|^2 \\
&\quad = e^{2\xi} - 2e^{\xi}(\cos\eta - \cos(\varphi+\eta)) + 2(1-\cos\varphi) \\
&\quad = e^{2\xi}(1 - 2(1-\cos\varphi)e^{-\xi}(\cos\eta - e^{-\xi}) - 2e^{-\xi}\sin\varphi\sin\eta) \\
&\quad \le \exp(2\xi - 2e^{-\xi}(1-\cos\varphi)(\cos\eta - e^{-\xi}) + 2e^{-\xi}|\sin\varphi||\sin\eta|),
\end{aligned}$$

which implies that

$$|1 + e^{\mathrm{i}\varphi}(e^z - 1)| \le \exp(\xi - (1-\cos\varphi)e^{-\xi}(\cos\eta - e^{-\xi}) + e^{-\xi}|\sin\varphi||\eta|).$$

This inequality (applied to $p_i z$) gives, by (3.5),

$$|P(ne^{\mathrm{i}\theta}, e^{\mathrm{i}\varphi})| \le \exp(\xi - (1-\cos\varphi)S_1 + |\sin\varphi|S_2), \tag{5.7}$$

where

$$\begin{aligned}
S_1 &:= \sum_j e^{-p_j\xi}(\cos(p_j\eta) - e^{-p_j\xi}), \\
S_2 &:= \sum_j e^{-p_j\xi}p_j|\eta| = |\tan\theta|u(\xi).
\end{aligned}$$

By (2.5), $u(\xi)^2 \le \xi v(\xi)$, and thus

$$|\sin\varphi|S_2 = O(|\theta\varphi|u(\xi)) = O(\sqrt{n\theta^2 v(\xi)\varphi^2}). \tag{5.8}$$

For $S_1$, let $c_8 > 0$ be chosen such that $\cos(x) - e^{-x} > 0$ on $(0, 2c_8]$ (e.g., $c_8 = 1/2$), and decompose the sum into three parts:

$$\begin{aligned}
S_1 &= \left(\sum_{p_j\xi \le c_8} + \sum_{c_8 < p_j\xi \le c_8/|\theta|} + \sum_{p_j\xi|\theta| > c_8}\right) e^{-p_j\xi}(\cos(p_j\eta) - e^{-p_j\xi}) \\
&=: T_1 + T_2 + T_3.
\end{aligned}$$

Consider first $T_1$. For each term in $T_1$ we have $p_j|\eta| \le p_j\xi \le c_8$, and since the function $x \mapsto (\cos x - e^{-x})/x$ extends to a continuous strictly positive function on $[0, c_8]$,

$$\cos(p_j\eta) - e^{-p_j\xi} \ge \cos(p_j\xi) - e^{-p_j\xi} \ge c_9 p_j\xi.$$

Furthermore, $e^{-p_j\xi} \geq e^{-c_8}$, and consequently

$$T_1 \geq c_10 \sum_{p_j\xi \leq c_8} p_j\xi.$$

For a term in $T_2$ we have either $p_j|\eta| \leq c_8$ and then

$$\cos(p_j\eta) - e^{-p_j\xi} \geq \cos(c_8) - e^{-c_8} > 0,$$

or $c_8 < p_j|\eta| = p_j\xi|\tan\theta| \leq 2p_j\xi|\theta| \leq 2c_8$ and then

$$\cos(p_j\eta) - e^{-p_j\xi} \geq \cos(p_j\eta) - e^{-p_j|\eta|} \geq c_{11}.$$

Consequently,

$$T_2 \geq c_{12} \sum_{c_8 < p_j\xi \leq c_8/|\theta|} e^{-p_j\xi}.$$

For $T_3$, which has at most $\xi|\theta|/c_8$ terms, we subtract $e^{-p_j\xi}$ from each term and use the trivial estimate $|\cos(p_j\eta) - e^{-p_j\xi} - 1| \leq 3$, finding

$$T_3 - \sum_{p_j\xi|\theta| > c_8} e^{-p_j\xi} = O(\xi|\theta|e^{-c_8/|\theta|}) = O(n\theta^4).$$

Combining these estimates, we obtain, using Lemma 4.1,

$$S_1 \geq c_{13}\Bigg(\sum_{p_j\xi \leq c_8} p_j\xi + \sum_{p_j\xi > c_8} e^{-p_j\xi}\Bigg) + O(n\theta^4) \geq c_{14}v(\xi) + O(n\theta^4). \tag{5.9}$$

The estimates (5.7), (5.8), (5.9), the Taylor expansion $\xi = n\cos\theta = n - n\theta^2/2 + O(n\theta^4)$, and the inequality $1 - \cos x \geq 2x^2/\pi^2$ for $x \in [-\pi, \pi]$ yield

$$|P(ne^{\mathrm{i}\theta}, e^{\mathrm{i}\varphi})| \leq \exp(n - n\theta^2/2 - c_{15}\varphi^2 v(\xi) + O(n\theta^4) + O(\sqrt{n\theta^2 v(\xi)\varphi^2})).$$

The required result (5.6) now follows, using $\sigma^2(n) \leq v(n) = O(\xi)$, provided $|\theta| \leq c_{16}$ and $n\theta^2 \leq c_{17}v(\xi)\varphi^2$. In both the remaining cases, the result follows from Proposition 5.3 if $c_6$ and $c_7$ are small enough. $\square$

5.2. *Local expansion for $P(ne^{\mathrm{i}\theta}, e^{\mathrm{i}\varphi})$.* We first rewrite (3.5) as

$$P(z, e^{\mathrm{i}\varphi}) = e^z \prod_j G(p_j z, \mathrm{i}\varphi), \tag{5.10}$$

where

$$G(z, \zeta) := 1 + (1 - e^{-z})(e^\zeta - 1) = e^\zeta - e^{\zeta - z} + e^{-z}, \tag{5.11}$$

with $z, \zeta \in \mathbb{C}$. We begin with an expansion of $G$.

LEMMA 5.7. *If* $|\arg z| \le \pi/3$ *and* $|\zeta| \le c_{18}$, *then*

$$G(z,\zeta) = \exp((1-e^{-z})\zeta + \tfrac{1}{2}(e^{-z}-e^{-2z})\zeta^2 + O(|\zeta|^3 e^{-\Re z}(1-e^{-\Re z}))).$$

PROOF. Choose $c_{18}$ such that $|e^{\zeta}-1| \le 1/4$ when $|\zeta| \le c_{18}$. Let

$$D := \{(z,\zeta) : |\arg z| \le \pi/3, |\zeta| \le c_{18}\}.$$

Then, for $(z,\zeta) \in D$,

$$|G(z,\zeta)-1| = |1-e^{-z}||e^{\zeta}-1| \le 2 \cdot \tfrac{1}{4} = \tfrac{1}{2}.$$

Hence, $g(z,\zeta) := \log G(z,\zeta)$ is well defined on $D$; moreover, $|G(z,\zeta)| \ge 1/2$ on $D$. Straightforward calculus yields, on $D$,

$$\frac{\partial}{\partial\zeta} g(z,\zeta) = \frac{(1-e^{-z})e^{\zeta}}{G(z,\zeta)} = 1 - \frac{e^{-z}}{G(z,\zeta)}, \tag{5.12}$$

$$\frac{\partial^2}{\partial\zeta^2} g(z,\zeta) = \frac{e^{-z}\frac{\partial}{\partial\zeta}G(z,\zeta)}{G(z,\zeta)^2} = \frac{e^{-z}(1-e^{-z})e^{\zeta}}{G(z,\zeta)^2}, \tag{5.13}$$

$$\frac{\partial^3}{\partial\zeta^3} g(z,\zeta) = O(|e^{-z}(1-e^{-z})|). \tag{5.14}$$

Since $|\arg z| \le \pi/3$, we see that if $|z| \le 1$, then $|1-e^{-z}| = O(|z|) = O(\Re z) = O(1-e^{-\Re z})$, and if $|z| \ge 1$, then $|1-e^{-z}| \le 2 = O(1-e^{-\Re z})$. Hence, in either case, $1-e^{-z} = O(1-e^{-\Re z})$, and by (5.14), we get $\frac{\partial^3}{\partial\zeta^3} g(z,\zeta) = O(e^{-\Re z}(1-e^{-\Re z}))$ on $D$. Moreover, $G(z,0)=1$ so $g(z,0)=0$, and, by (5.12) and (5.13), $\frac{\partial}{\partial\zeta} g(z,0) = 1-e^{-z}$ and $\frac{\partial^2}{\partial\zeta^2} g(z,0) = e^{-z}(1-e^{-z})$. This proves the lemma. □

PROPOSITION 5.8. *If* $|\theta| \le \pi/3$ *and* $|\varphi| \le c_{18}$, *then*

$$P(ne^{\mathrm{i}\theta}, e^{\mathrm{i}\varphi}) = \exp(ne^{\mathrm{i}\theta} + \mu(n)\mathrm{i}\varphi - u(n)\varphi\theta - \tfrac{1}{2}v(n)\varphi^2 + O(v(n\cos\theta)|\varphi|^3 + n\theta^2|\varphi|)).$$

PROOF. Let $z := ne^{\mathrm{i}\theta}$ and $\xi := \Re z = n\cos\theta$. By (5.10) and Lemma 5.7,

$$P(z, e^{\mathrm{i}\varphi}) = \exp(z + \mu(z)\mathrm{i}\varphi - \tfrac{1}{2}v(z)\varphi^2 + O(v(\xi)|\varphi|^3)).$$

Observe that $\mu'(z) = \sum_j p_j e^{-p_j z} = u(z)/z$ and $\mu''(z) = -\sum_j p_j^2 e^{-p_j z}$. Thus, by a Taylor expansion and Lemmas 5.4 and 5.5,

$$\mu(z) = \mu(n) + u(n)\mathrm{i}\theta + O(n\theta^2). \tag{5.15}$$

Similarly, using the inequality $u^2(x)/x \le v(x)$,

$$|v'(z)| = \left|\sum_j (2p_j e^{-2p_j z} - p_j e^{-p_j z})\right| \le 3u(\xi)/\xi \le 3(v(\xi)/\xi)^{1/2},$$

and thus

$$v(z) = v(n) + O(|\theta|(nv(\xi))^{1/2}).$$

The desired result follows from these estimates and the inequality $(nv(\xi))^{1/2} \times |\theta|\varphi^2 \le v(\xi)|\varphi|^3 + n\theta^2|\varphi|$. $\square$

5.3. *Proof of Theorem 2.4 when* $\sum_{p_j n>1} p_j \ge 1/2$. The remaining analysis is straightforward. We assume that $\sigma^2(n) \ge 1$, since otherwise the result is trivial. Recall that $\theta_0 := n^{-1/2}\sigma(n)^{1/3} \le n^{-1/3}$ and $\varphi_0 := \sigma(n)^{-2/3}$. We assume that $\sum_{p_j n>1} p_j \ge 1/2$, and thus $\sigma^2(n) \asymp v(n)$ by Proposition 4.3.

We start from (5.1) and split the integral into the three parts in (5.2):

$$\begin{aligned}\mathbb{P}(Z_n = m) \\ &= \frac{n!n^{-n}}{(2\pi)^2}\left(\iint_{\substack{|\theta|\le\theta_0\\|\varphi|\le\varphi_0}} + \iint_{\substack{|\theta|\le\theta_0\\ \varphi_0<|\varphi|\le\pi}} + \iint_{\substack{\theta_0<|\theta|\le\pi\\|\varphi|\le\pi}}\right)e^{-\mathrm{i}m\varphi-\mathrm{i}n\theta} \\ &\qquad\qquad \times P(ne^{\mathrm{i}\theta}, e^{\mathrm{i}\varphi})\,d\theta\,d\varphi \\ &=: J_1 + J_2 + J_3.\end{aligned}$$

Observe first that, by Stirling's formula,

$$\frac{n!n^{-n}}{(2\pi)^2} = (2\pi)^{-3/2}\sqrt{n}e^{-n}(1 + O(1/n)) \tag{5.16}$$

and thus

$$\frac{n!n^{-n}}{(2\pi)^2} = O(\sqrt{n}e^{-n}). \tag{5.17}$$

Obviously, by Proposition 5.3 and (5.17),

$$\begin{aligned}J_3 &= O\left(\sqrt{n}\int_{\theta_0}^{\infty} e^{-c_5 n\theta^2}\,d\theta\right) = O(n^{-1/2}\theta_0^{-1}e^{-c_5 n\theta_0^2}) \\ &= O(e^{-c_5\sigma(n)^{2/3}}) = O(\sigma^{-2}(n)).\end{aligned} \tag{5.18}$$

On the other hand, Proposition 5.6 gives, for $n \ge C_3$ and $|\theta| \le \theta_0$,

$$\begin{aligned}J_2 &= O\left(\sqrt{n}\int_0^{\infty} e^{-c_6 n\theta^2}\,d\theta \int_{\varphi_0}^{\infty} e^{-c_7\sigma^2(n)\varphi^2}\,d\varphi\right) \\ &= O(e^{-c_7\sigma^2(n)\varphi_0^2}) = O(e^{-c_7\sigma(n)^{2/3}}) = O(\sigma^{-2}(n)).\end{aligned} \tag{5.19}$$

We turn to $J_1$, the main term. If $n \ge C_4$ and $\sigma^2(n) \ge C_5$, then Proposition 5.8 applies when $|\theta| \le \theta_0$ and $|\varphi| \le \varphi_0$, and shows, together with Proposition 4.3 and Lemma 5.5(iv), that

$$P(ne^{\mathrm{i}\theta}, e^{\mathrm{i}\varphi}) = \exp(n + \mathrm{i}n\theta - \tfrac{1}{2}n\theta^2 + \mu(n)\mathrm{i}\varphi - u(n)\varphi\theta - \tfrac{1}{2}v(n)\varphi^2 + R(\theta,\varphi)),$$

where

$$R(\theta,\varphi)=O(\sigma^2(n)|\varphi|^3+n\theta^2|\varphi|+n|\theta|^3)=O(1).$$

Let

$$K(\theta,\varphi):=\exp(\mathrm{i}(\mu(n)-m)\varphi-\tfrac12 n\theta^2-u(n)\theta\varphi-\tfrac12 v(n)\varphi^2). \tag{5.20}$$

Then, by the inequality $|e^z-1|\le |z|e^{|z|}$,

$$\begin{aligned} J_1&=\frac{n!n^{-n}e^n}{(2\pi)^2}\int_{-\varphi_0}^{\varphi_0}\int_{-\theta_0}^{\theta_0}K(\theta,\varphi)e^{R(\theta,\varphi)}\,d\theta\,d\varphi\\ &=\frac{n!n^{-n}e^n}{(2\pi)^2}\int_{-\varphi_0}^{\varphi_0}\int_{-\theta_0}^{\theta_0}K(\theta,\varphi)(1+O(|R(\theta,\varphi)|)\,d\theta\,d\varphi. \end{aligned}$$

We first estimate the error term. Observe that

$$|K(\theta,\varphi)|=\exp(-\tfrac12 n\theta^2-u(n)\theta\varphi-\tfrac12 v(n)\varphi^2)=\exp(-\tfrac12 A(\sqrt{n}\theta,\sqrt{v(n)}\varphi)),$$

where $A$ is the quadratic form

$$A(x,y):=x^2+y^2+2\frac{u(n)}{\sqrt{nv(n)}}xy.$$

Since

$$\left(\frac{u(n)}{\sqrt{nv(n)}}\right)^2=\frac{u(n)^2/n}{v(n)}=\frac{v(n)-\sigma^2(n)}{v(n)}\le 1-c_{19}$$

by Proposition 4.3(ii), we see that $A(x,y)\ge c_{20}(x^2+y^2)$, implying that

$$|K(\theta,\varphi)|\le e^{-c_{21}n\theta^2-c_{21}\sigma^2(n)\varphi^2}. \tag{5.21}$$

It follows that

$$\begin{aligned} &\frac{n!n^{-n}e^n}{(2\pi)^2}\int_{-\varphi_0}^{\varphi_0}\int_{-\theta_0}^{\theta_0}K(\theta,\varphi)|R(\theta,\varphi)|\,d\theta\,d\varphi\\ &\qquad\le C_6\sqrt{n}\int_{-\infty}^{\infty}\int_{-\infty}^{\infty}(\sigma^2(n)|\varphi|^3+n\theta^2|\varphi|+n|\theta|^3)e^{-c_{21}n\theta^2-c_{21}\sigma^2(n)\varphi^2}\,d\theta\,d\varphi\\ &\qquad=O(\sigma(n)^{-2}+n^{-1/2}\sigma(n)^{-1})\\ &\qquad=O(\sigma(n)^{-2}). \end{aligned}$$

It remains only to evaluate the integral of $K$ over the remaining region. The estimate (5.21) implies, arguing as for $J_3$ and $J_2$ in (5.18) and (5.19), that

$$\frac{n!n^{-n}e^n}{(2\pi)^2}\left(\iint_{\substack{|\varphi|>\varphi_0\\|\theta|\le\theta_0}}+\iint_{\substack{-\infty<\varphi<\infty\\|\theta|>\theta_0}}\right)|K(\theta,\varphi)|\,d\theta\,d\varphi=O(\sigma(n)^{-2}).$$

Collecting the estimates above, we get

$$\mathbb{P}(Z_n = m) = \frac{n!n^{-n}e^n}{(2\pi)^2}\int_{-\infty}^{\infty}\int_{-\infty}^{\infty} K(\theta,\varphi)\,d\theta\,d\varphi + O(\sigma(n)^{-2}). \tag{5.22}$$

Since (5.20) can be rewritten as

$$K(\theta,\varphi) = \exp(-\mathrm{i}(m-\mu(n))\varphi - \tfrac{1}{2}n(\theta + (u(n)/n)\varphi)^2 - \tfrac{1}{2}\sigma^2(n)\varphi^2),$$

it then follows that

$$\frac{1}{2\pi}\int_{-\infty}^{\infty}\int_{-\infty}^{\infty} K(\theta,\varphi)\,d\theta\,d\varphi = n^{-1/2}\sigma(n)^{-1}\exp\Big(-\frac{(m-\mu(n))^2}{2\sigma^2(n)}\Big),$$

which, together with (5.22) and (5.16), completes the proof; note that if $m = \lfloor \mu(n) + x\sigma(n)\rfloor$, then $(m-\mu(n))/\sigma(n) = x + O(1/\sigma(n))$, and that $x \mapsto e^{-x^2}/2$ has bounded derivatives. The assumptions above that $n$ and $\sigma^2(n)$ be large are harmless since $n \geq \sigma^2(n)$ and the result is trivial for $\sigma^2(n) \leq C_7$, for any fixed $C_7$.

REMARK 5.1 (*Central limit theorem*). If one is interested in proving only the central limit theorem, then Propositions 5.3 and 5.8 suffice. If, moreover, a Berry–Esseen bound is desired, then Proposition 5.6 is needed for $|\varphi| \leq \varepsilon$ for some $\varepsilon > 0$.

## 6. Proof of Theorem 2.4 when $\sum_{p_j n \leq 1} p_j \geq 1/2$.

We consider in this section the case when $\sum_{p_j n\leq 1} p_j \geq 1/2$. Our underlying idea is then to study $n - Z_n$ instead of $Z_n$, and the corresponding Poissonization $\tilde{Z}(n)$; see Remark 3.2 and recall that $\mathrm{Var}(\tilde{Z}(n)) = \tilde{v}(n) \asymp \sigma^2(n)$. We find $\mathbb{P}(Z_n = m) = \mathbb{P}(n - Z_n = n - m)$ by extracting coefficients in $P(e^{\mathrm{i}\varphi}\lambda, e^{-\mathrm{i}\varphi})$; see (3.6). This yields the integral formula

$$\begin{aligned}\mathbb{P}&(n - Z_n = n-m)\\ &= \frac{n!n^{-n}}{(2\pi)^2}\int_{-\pi}^{\pi}\int_{-\pi}^{\pi} e^{-(n-m)\mathrm{i}\varphi - n\mathrm{i}\theta} P(ne^{\mathrm{i}\theta+\mathrm{i}\varphi}, e^{-\mathrm{i}\varphi})\,d\theta\,d\varphi.\end{aligned} \tag{6.1}$$

Note that this formula also follows directly from (5.1) by a simple change of variables; there is thus no formal need of $\tilde{Z}$ and the motivation above.

The analysis of this double integral is very similar to the one in Section 5; the main difference is that the occurrences of $v(n)$ in our estimates have to be replaced by $\tilde{v}(n)$.

### 6.1. *Estimates for* $|P(ne^{\mathrm{i}\theta+\mathrm{i}\varphi}, e^{-\mathrm{i}\varphi})|$.

We begin with a companion to Lemma 5.2.

LEMMA 6.1. *If* $0 \leq r \leq 1$ *and* $|\theta| \leq \pi$, *then*

$$1 + |e^{re^{\mathrm{i}\theta}} - 1| \leq e^{r - c_{22}r^2\theta^2}.$$

Proof. By Lemma 5.1, we have (5.4) and thus

$$\begin{aligned} e^{-r}(1+|e^{re^{i\theta}}-1|) &\le e^{-r}+(1-e^{-r})e^{-c_4 r\theta^2} \\ &= 1-(1-e^{-r})(1-e^{-c_4 r\theta^2}) \\ &\le \exp(-(1-e^{-r})(1-e^{-c_4 r\theta^2})), \end{aligned}$$

and the result follows. □

Lemma 6.2. *Uniformly for* $|\theta| \le \pi$ *and* $-\infty < \varphi < \infty$,

$$|P(ne^{i\theta}, e^{i\varphi})| \le e^{n - c_{23}\tilde{v}(n)\theta^2}.$$

Proof. A simple consequence of (3.5), Lemmas 5.2 and 6.1 and (4.5). □

Lemma 6.3. (i) *If* $r \ge 0$, $|\theta| \le \pi$ *and* $|\zeta| = 1$, *then*

$$|\zeta + e^{\zeta r e^{i\theta}} - 1| \le e^r.$$

(ii) *If, furthermore,* $0 \le r \le 1$, *then*

$$|\zeta + e^{\zeta r e^{i\theta}} - 1| \le e^{r - c_{24} r\theta^2}.$$

Proof. Expanding the function $\zeta + e^{\zeta r e^{i\theta}} - 1$ at $r = 0$ gives

$$\begin{aligned} |\zeta + e^{\zeta r e^{i\theta}} - 1| &\le |\zeta + \zeta r e^{i\theta}| + \sum_{k\ge 2} \frac{|\zeta r e^{i\theta}|^k}{k!} \\ &= |1 + re^{i\theta}| + \sum_{k\ge 2}\frac{r^k}{k!} \qquad (6.2) \\ &= |1 + re^{i\theta}| + e^r - r - 1. \end{aligned}$$

Part (i) follows immediately. On the other hand, since

$$\begin{aligned} |1+re^{i\theta}|^2 &= (1+r)^2 - 2r(1-\cos\theta) \\ &\le (1+r)^2 - 4c_4 r\theta^2, \end{aligned}$$

we have the inequality

$$|1+re^{i\theta}| \le 1 + r - c_{25} r\theta^2 \qquad (r \in [0,1]).$$

This together with (6.2) yields

$$|\zeta + e^{\zeta r e^{i\theta}} - 1| \le e^r - c_{25} r\theta^2 \le e^r(1 - c_{24} r\theta^2) \le e^{r - c_{24} r\theta^2},$$

uniformly for $r \in [0,1]$. □

The next proposition is the analogue of Proposition 5.3 when $\sum_{p_j n \le 1} p_j \ge 1/2$.

PROPOSITION 6.4. *Assume that* $\sum_{p_j n \le 1} p_j \ge 1/2$. *Then uniformly for* $|\theta| \le \pi$ *and* $-\infty < \varphi < \infty$,

$$|P(ne^{\mathrm{i}(\theta+\varphi)}, e^{-\mathrm{i}\varphi})| \le e^{n - c_{26} n\theta^2}.$$

PROOF. By (3.5) and Lemma 6.3 with $\zeta = e^{\mathrm{i}\varphi}$,

$$\begin{aligned} |P(ne^{\mathrm{i}(\theta+\varphi)}, e^{-\mathrm{i}\varphi})| &= \prod_j |e^{\mathrm{i}\varphi} + e^{p_j n e^{\mathrm{i}(\theta+\varphi)}} - 1| \\ &\le \biggl(\prod_{p_j n \le 1} e^{p_j n - c_{24} p_j n \theta^2}\biggr)\biggl(\prod_{p_j n > 1} e^{p_j n}\biggr) \\ &= \exp\biggl(n - c_{24}\theta^2 \sum_{p_j n \le 1} p_j n\biggr). \end{aligned}$$ □

The corresponding analogue of Proposition 5.6 is the following.

PROPOSITION 6.5. *Assume that* $\sum_{p_j n \le 1} p_j \ge 1/2$. *Then uniformly for* $|\theta| \le \pi$ *and* $|\varphi| \le \pi$,

$$|P(ne^{\mathrm{i}(\theta+\varphi)}, e^{-\mathrm{i}\varphi})| \le \exp(n - c_{27}(n\theta^2 + \sigma^2(n)\varphi^2)).$$

PROOF. If $|\varphi| \le 2|\theta|$, then $n\theta^2 + \sigma^2(n)\varphi^2 = O(n\theta^2)$, and the result follows by Proposition 6.4.

On the other hand, if $|\varphi| \ge 2|\theta|$, then $|\theta + \varphi| \le \frac{3}{2}|\varphi| \le \frac{3}{2}\pi$. Note that Lemma 6.2 extends to $|\theta| \le \frac{3}{2}\pi$ (with a new $c_{23}$) since if $\pi < |\theta| \le \frac{3}{2}\pi$, we may replace $\theta$ by $\theta \pm 2\pi$. Hence, by Lemma 6.2 and Proposition 6.4,

$$e^{-n}|P(ne^{\mathrm{i}(\theta+\varphi)}, e^{-\mathrm{i}\varphi})| \le \exp(-\tfrac{1}{2}(c_{23}\tilde{v}(n)(\theta+\varphi)^2 + c_{26} n\theta^2)),$$

and the result follows because $\sigma^2(n) = O(\tilde{v}(n))$, $\varphi^2 \le 2\theta^2 + 2(\theta+\varphi)^2$, and thus $n\theta^2 + \sigma^2(n)\varphi^2 = O(n\theta^2 + \tilde{v}(n)(\theta+\varphi)^2)$. □

6.2. *Local expansion for* $P(ne^{\mathrm{i}\theta+\mathrm{i}\varphi}, e^{-\mathrm{i}\varphi})$. We turn now to a local expansion of $P(ne^{\mathrm{i}\theta+\mathrm{i}\varphi}, e^{-\mathrm{i}\varphi})$. We will use

$$P(ze^{\mathrm{i}\varphi}, e^{-\mathrm{i}\varphi}) = \prod_j (1 + e^{p_j z e^{\mathrm{i}\varphi} - \mathrm{i}\varphi} - e^{-\mathrm{i}\varphi}) = e^z \prod_j H(p_j z, \mathrm{i}\varphi), \tag{6.3}$$

where we define

$$H(z,\zeta) = e^{-z}(1 + e^{ze^{\zeta} - \zeta} - e^{-\zeta}) = e^{z(e^{\zeta}-1)-\zeta} + e^{-z}(1 - e^{-\zeta}). \tag{6.4}$$

LEMMA 6.6. *If $|z| \le 1$ and $|\zeta| \le c_{28}$, then*

$$H(z,\zeta) = \exp((z-1+e^{-z})\zeta + \tfrac{1}{2}(z+e^{-z}-e^{-2z}-2ze^{-z})\zeta^2 + O(|z^2\zeta^3|)).$$

PROOF. Note first that $H(z,0)=1$. Hence, for $|z|\le 1$ and $|\zeta|\le c_{28}$, we have $|H(z,0)-1|\le 1/2$. Thus $h(z,\zeta) := \log H(z,\zeta)$ is well defined in the domain $D := \{(z,\zeta) : |z| \le 1, |\zeta| \le c_{28}\}$, with all its derivatives bounded and $h(z,0)=0$. Moreover, $h(0,\zeta)=0$ because $H(0,\zeta)=1$. Also, by

$$\frac{\partial}{\partial z} H(z,\zeta) = (e^{\zeta}-1)e^{z(e^{\zeta}-1)-\zeta} - e^{-z}(1-e^{-\zeta}),$$

we have $\frac{\partial}{\partial z} h(0,\zeta) = \frac{\partial}{\partial z} H(0,\zeta) = 0$.

Consequently,

$$\frac{\partial^3}{\partial\zeta^3} h(0,\zeta) = 0 \quad \text{and} \quad \frac{\partial}{\partial z}\frac{\partial^3}{\partial\zeta^3} h(0,\zeta) = \frac{\partial^3}{\partial\zeta^3}\frac{\partial}{\partial z} h(0,\zeta) = 0,$$

and a Taylor expansion in $z$ yields, for $(z,\zeta)\in D$,

$$\frac{\partial^3}{\partial\zeta^3} h(z,\zeta) = O(|z|^2).$$

Hence, by another Taylor expansion, now in $\zeta$, for $(z,\zeta)\in D$,

$$h(z,\zeta) = \frac{\partial}{\partial\zeta} h(z,0)\zeta + \frac{1}{2}\frac{\partial^2}{\partial\zeta^2} h(z,0)\zeta^2 + O(|z^2\zeta^3|),$$

and the result follows, with the values of $\frac{\partial}{\partial\zeta} h(z,0)$ and $\frac{\partial^2}{\partial\zeta^2} h(z,0)$ obtained by straightforward calculus. □

LEMMA 6.7. *If $|\arg z| \le \pi/4$ and $|\zeta| \le c_{29}$, then*

$$H(z,\zeta) = \exp((z-1+e^{-z})\zeta + \tfrac{1}{2}(z+e^{-z}-e^{-2z}-2ze^{-z})\zeta^2 + O(|z\zeta^3|)).$$

PROOF. By the definitions (6.4) and (5.11), with $w := ze^{\zeta}$,

$$H(z,\zeta) = e^{-z+w}(e^{-\zeta}+e^{-w}-e^{-w-\zeta}) = e^{-z+w}G(w,-\zeta).$$

Thus, by Lemma 5.7, for $|\arg z|\le \pi/4$ and $|\zeta|\le c_{29}$,

$$\begin{aligned} H(z,\zeta) = \exp(z(e^{\zeta}-1) - (1-e^{-w})\zeta \\ + \tfrac{1}{2}(e^{-w}-e^{-2w})\zeta^2 + O(|\zeta|^3)). \end{aligned} \tag{6.5}$$

Moreover, for $|\arg z|\le\pi/4$ and $|\zeta|\le c_{29}$,

$$\begin{aligned} \frac{\partial}{\partial\zeta} e^{-ze^{\zeta}} &= -ze^{\zeta}e^{-ze^{\zeta}}, \\ \frac{\partial^2}{\partial\zeta^2} e^{-ze^{\zeta}} &= ((ze^{\zeta})^2 - ze^{\zeta})e^{-ze^{\zeta}} = O((|z|^2+|z|)e^{-c_{30}|z|}) = O(1), \end{aligned}$$

and thus

$$e^{-w} = e^{-z} - ze^{-z}\zeta + O(|\zeta|^2).$$

The result follows by substituting this and $e^{\zeta} - 1 = \zeta + \frac{1}{2}\zeta^2 + O(|\zeta|^3)$ in (6.5), provided $|z| \geq 1$. The case $|z| < 1$ is a consequence of Lemma 6.6. □

LEMMA 6.8. *Uniformly for* $\Re z \geq 1$ *in the sector* $|\arg z| \leq \pi/4$,

$$|\tilde{v}'(z)| = O(\tilde{v}(|z|)/|z|).$$

PROOF. Let $\tau(z) := z + e^{-z}(1 - e^{-z} - 2z)$. Then $\tilde{v}(z) = \sum_j \tau(p_j z)$ and $\tilde{v}'(z) = \sum_j p_j \tau'(p_j z)$. Since $\tau'(0) = 0$, we see that $\tau'(z) = O(|z|)$ for $|z| \leq 1$. Furthermore, it is easily seen that $\tau'(z) = O(1)$ when $|z| \geq 1$ in the sector $|\arg z| \leq \pi/3$. Hence, using (4.5),

$$\begin{aligned}|\tilde{v}'(z)| &\leq C_8 \sum_j p_j (p_j|z| \wedge 1) \\ &= C_8 x |z|^{-1} \sum_j (p_j|z|)^2 \wedge (p_j|z|) \leq C_9 |z|^{-1} \tilde{v}(|z|).\end{aligned}$$

This completes the proof. □

The next result gives the analogue of Proposition 5.8.

PROPOSITION 6.9. *If* $|\theta| \leq \pi/4$ *and* $|\varphi| \leq c_{31}$, *then*

$$\begin{aligned}P(ne^{\mathrm{i}\theta+\mathrm{i}\varphi}, e^{-\mathrm{i}\varphi}) = \exp(&ne^{\mathrm{i}\theta} + (n - \mu(n))\mathrm{i}\varphi - (n - u(n))\varphi\theta - \tfrac{1}{2}\tilde{v}(n)\varphi^2 \\ &+ O(\tilde{v}(n)|\varphi|^3 + n\theta^2|\varphi|)).\end{aligned} \tag{6.6}$$

PROOF. Let $z := ne^{\mathrm{i}\theta}$. It follows from (6.3), Lemmas 6.6 and 6.7, and (2.2)–(2.6) together with (4.5) that, assuming $|\varphi| \leq c_{31}$,

$$P(ze^{\mathrm{i}\varphi}, e^{-\mathrm{i}\varphi}) = \exp(z + \mathrm{i}(z - \mu(z))\varphi - \tfrac{1}{2}\tilde{v}(z)\varphi^2 + O(\tilde{v}(n)|\varphi|^3)).$$

By (5.15),

$$z - \mu(z) = n - \mu(n) + \mathrm{i}\theta(n - u(n)) + O(n\theta^2).$$

On the other hand, by Lemma 6.8, we also have

$$\tilde{v}(z) = \tilde{v}(n) + O(|\theta|\tilde{v}(n)),$$

and the result (6.6) follows, in view of the inequalities $\tilde{v}(n)|\theta|\varphi^2 \leq \tilde{v}(n)|\varphi|^3 + \tilde{v}(n)\theta^2|\varphi|$ and $\tilde{v}(n) \leq n$. □

6.3. *Proof of Theorem* 2.4 *when* $\sum_{p_jn\le 1} p_j \ge 1/2$. The analysis of (6.1) is essentially the same as was done for (5.1) in Section 5.3, now using Propositions 6.5 and 6.9 and the relation

$$\tilde{v}(n) = (n - u(n))^2/n + \sigma^2(n),$$

which follows from (2.5) and (2.6). We omit the details.

## 7. Proofs of Theorem 2.1 and Corollary 2.5.

Proof of Theorem 2.1. We may assume that $\sigma_F^2(n) \ge 1$ and $\operatorname{Var}(Z_{n,F}) \ge 1$, since otherwise $\operatorname{Var}(Z_{n,F}) \le C_{10}$ by Theorem 2.3 and the result is trivial. Then (2.1) is a simple consequence of Theorems 2.4 and 2.3. □

Proof of Corollary 2.5. (i) $\Longleftrightarrow$ (ii). An immediate consequence of Theorem 2.3.

(i) $\Longrightarrow$ (iii). By Theorem 2.1.

(ii) $\Longrightarrow$ (iv). By Theorem 2.4.

(iii) $\Longrightarrow$ (v) and (iv) $\Longrightarrow$ (v). Trivial.

(v) $\Longrightarrow$ (i). (This part is standard and uses the fact that $Z_\nu$ assumes only integer values.) If (v) holds, let $Z'_\nu$ be an independent copy of $Z_\nu$. Then

$$(Z_\nu - Z'_\nu)/\beta_\nu \xrightarrow{\mathrm{d}} N(0,2). \tag{7.1}$$

If (i) fails, then there is a subsequence $(n_\nu, F_\nu)_{\nu\in N'}$, along which $\sigma^2_{n_\nu,F_\nu}$ is bounded; we consider that subsequence only, and let $B := \sup_{\nu\in N'} \beta_\nu$.

If $B = \infty$, there is a subsubsequence along which $\beta_\nu \to \infty$, but this implies $\mathbb{E}((Z_\nu - Z'_\nu)/\beta_\nu)^2 = 2\sigma^2_{n_\nu,F_\nu}/\beta_\nu^2 \to 0$, and thus $(Z_\nu - Z'_\nu)/\beta_\nu \xrightarrow{\mathrm{d}} 0$ along the subsubsequence, which contradicts (7.1).

On the other hand, if $B < \infty$, then $\mathbb{P}((Z_\nu - Z'_\nu)/\beta_\nu \in [1/4B, 1/2B]) = 0$ for all $\nu \in N'$ since $Z_\nu - Z'_\nu$ is integer-valued, which again contradicts (7.1). □

**8. Limit laws when the variance is small or bounded.** We briefly consider the possible limit laws for a sequence of random variables $Z_n = Z_{n,F_n}$ with bounded variances. [Recall that Corollary 2.5 shows that $Z_n$ is asymptotically normal in the opposite case when $\operatorname{Var}(Z_n) \to \infty$.] By Theorem 2.3, this assumption is equivalent to $\sigma^2(n) = O(1)$, and according to Proposition 4.3 and Remark 4.1, we consider the following two cases:

(i) $\sum_{p_jn\le 1} p_j \ge 1/2$, $\tilde{v}(n) = O(1)$, $v(n) \asymp n$;

(ii) $\sum_{p_jn>1} p_j \ge 1/2$, $v(n) = O(1)$, $\tilde{v}(n) \asymp n$.

In both cases we can use the same Poissonization procedure as above, the proofs being indeed much simpler. However, for more methodological interests, we use the coupling argument mentioned in Remark 3.1. We say that an event holds whp (*with high probability*), if it holds with probability tending to 1 as $n \to \infty$.

PROPOSITION 8.1. (i) *If $n \to \infty$ with $\tilde{v}(n) = O(1)$, then whp $n - Z_n = \tilde{Z}(n)$.*

(ii) *If $n \to \infty$ with $v(n) = O(1)$, then whp $Z_n = Z(n)$.*

PROOF. Let $\lambda_\pm := n \pm n^{2/3}$. Then, whp, $N(\lambda_-) \le n \le N(\lambda_+)$, and thus $Z(\lambda_-) \le Z_n \le Z(\lambda_+)$ and $\tilde{Z}(\lambda_-) \le n - Z_n \le \tilde{Z}(\lambda_+)$. Moreover, $Z(\lambda_-) \le Z(n) \le Z(\lambda_+)$ and $\tilde{Z}(\lambda_-) \le \tilde{Z}(n) \le \tilde{Z}(\lambda_+)$. Consequently, it suffices to show that whp $Z(\lambda_-) = Z(\lambda_+)$ in case (ii) and $\tilde{Z}(\lambda_-) = \tilde{Z}(\lambda_+)$ in case (i).

In case (i) we have for $\lambda_- \le \lambda \le \lambda_+$, using (4.5),

$$\begin{aligned}\frac{d}{d\lambda}\mathbb{E}(\tilde{Z}(\lambda)) &= \frac{d}{d\lambda}(\lambda - \mu(\lambda)) = 1 - \sum_j p_j e^{-p_j\lambda} = \sum_j p_j(1 - e^{-p_j\lambda}) \\ &\le \sum_j p_j(p_j\lambda \wedge 1) = O\left(\sum_j (p_j^2 n \wedge p_j)\right) = O(\tilde{v}(n)/n),\end{aligned}$$

and thus

$$\begin{aligned}\mathbb{P}(\tilde{Z}(\lambda_+) \ne \tilde{Z}(\lambda_-)) &\le \mathbb{E}(\tilde{Z}(\lambda_+) - \tilde{Z}(\lambda_-)) = O((\lambda_+ - \lambda_-)\tilde{v}(n)/n) \\ &= O(n^{-1/3}\tilde{v}(n)) = o(1).\end{aligned}$$

In case (ii) we have by Lemma 5.4, for $\lambda \ge 2$,

$$\frac{d}{d\lambda}\mathbb{E}(Z(\lambda)) = \sum_j p_j e^{-p_j\lambda} = \frac{u(\lambda)}{\lambda} = O(\lambda^{-1}\log(\lambda)v(\lambda)).$$

By Lemma 5.5 we thus have, for $\lambda \in [\lambda_-, \lambda_+]$ (and $n \ge 2$)

$$\frac{d}{d\lambda}\mathbb{E}(Z(\lambda)) = O(n^{-1}\log(n)v(n))$$

and accordingly

$$\mathbb{P}(Z(\lambda_+) \ne Z(\lambda_-)) \le \mathbb{E}(Z(\lambda_+) - Z(\lambda_-)) = O(n^{-1/3}\log(n)v(n)) = o(1). \quad \square$$

Limit results can now be obtained from the representations $Z(n) = \sum_j \mathbf{1}_{\{U_j \ge 1\}}$ and $\tilde{Z}(n) = \sum_j \tilde{U}_j$ with independent summands given in Section 3.

We consider in detail two simple cases leading to Poisson limit laws. Both cases are marked by the property that there are no $p_j$ of order $1/n$; compare Chistyakov [6] and Kolchin, Sevast'yanov and Chistyakov [23], III.3.

THEOREM 8.2. *Suppose that* $\frac{1}{2}\sum_j (p_j n)^2 \to \lambda < \infty$ *and that* $\max_j p_j n \to 0$*; then* $n - Z_n \xrightarrow{\mathrm{d}} \mathrm{Po}(\lambda)$.

PROOF. We have, by (4.5), $\tilde v(n) = O(\sum_j (p_j n)^2) = O(1)$, so Proposition 8.1(i) applies. Further, we have by Section 3, with $U_j \sim \mathrm{Po}(p_j n)$ independent,

$$\sum_j \mathbb{P}(\tilde U_j \ge 2) = \sum_j \mathbb{P}(U_j \ge 3) \le \sum_j (p_j n)^3 \le \max(p_j n) \sum_j (p_j n)^2 \to 0.$$

Thus, whp,

$$n - Z_n = \tilde Z(n) = \sum_j \tilde U_j = \sum_j \tilde I_j,$$

where $\tilde I_j := \mathbf{1}_{\{\tilde U_j = 1\}} = \mathbf{1}_{\{U_j = 2\}} \sim \mathrm{Be}(\frac{1}{2}(p_j n)^2 e^{-p_j n})$ are independent Bernoulli distributed variables. We have, as $n \to \infty$, $\max_j \mathbb{E}(\tilde I_j) \to 0$ and

$$\sum_j \mathbb{E}(\tilde I_j) = \sum_j \tfrac{1}{2}(p_j n)^2 e^{-p_j n} = \sum_j \tfrac{1}{2}(p_j n)^2 + O\left(\sum_j (p_j n)^3\right) \to \lambda.$$

Hence $\sum_j \tilde I_j \xrightarrow{\mathrm{d}} \mathrm{Po}(\lambda)$ by a standard result; see [14, 25] or, for example, [3], Theorem 2.M. □

THEOREM 8.3. *Suppose that:*

(i) $\sum_{p_j n \le 1} p_j n \to \lambda_1 \in [0, \infty)$,
(ii) $\sum_{p_j n > 1} e^{-p_j n} \to \lambda_2 \in [0, \infty)$,
(iii) $\sup_j (p_j n \wedge (p_j n)^{-1}) \to 0$,

*and let* $m := \#\{j : p_j > 1/n\}$. *Then* $Z_n - m \xrightarrow{\mathrm{d}} W_1 - W_2$, *where* $W_i \sim \mathrm{Po}(\lambda_i)$ *are independent.*

*Note that* $m$ *depends on* $n$ *and the* $p_j$*'s.*

PROOF. We have, by (4.4), $v(n) = O(1)$, so Proposition 8.1(ii) applies and whp $Z_n = Z(n) = \sum_j U_j$, where, by Section 3, $U_j \sim \mathrm{Be}(1 - \exp(-p_j n))$ are independent. Hence, whp,

$$Z_n - m = \sum_{p_j n \le 1} U_j - \sum_{p_j n > 1} (1 - U_j),$$

where the two sums of independent Bernoulli variables are independent. We have, as $n \to \infty$,

$$\sup_{p_j n \le 1} \mathbb{E}(U_j) \le \sup_{p_j n \le 1} p_j n \le \sup_j (p_j n \wedge (p_j n)^{-1}) \to 0$$

and

$$\sum_{p_jn\le 1}\mathbb{E}(U_j)=\sum_{p_jn\le 1}(1-e^{-p_jn})=\sum_{p_jn\le 1}p_jn+O\Bigg(\sum_{p_jn\le 1}(p_jn)^2\Bigg)\to\lambda_1,$$

because

$$\sum_{p_jn\le 1}(p_jn)^2\le\sup_j(p_jn\wedge(p_jn)^{-1})\sum_j p_jn\to 0.$$

Hence $\sum_{p_jn\le 1}U_j\xrightarrow{\mathrm{d}}W_1\sim\mathrm{Po}(\lambda_1)$, again by [14, 25] or, for example, [3], Theorem 2.M.

Similarly, $\sum_{p_jn>1}(1-U_j)\xrightarrow{\mathrm{d}}W_2\sim\mathrm{Po}(\lambda_2)$. $\square$

REMARK 8.1 (*Poisson approximation*). We can derive more precise local limit theorems by modifying our proof for Theorem 2.4; the proof is indeed much simpler and omitted here.

Theorems 8.2 and 8.3 extend to the general case when some $p_j$ is of the order $1/n$, but the limit distributions become more complicated. Consider first case (i), with $\tilde{v}(n)=O(1)$. We may assume that, for each $n$, $p_1\ge p_2\ge\cdots$; by (4.5), $p_1n=O(1)$ and we may by taking a subsequence [of $(n,F)$] assume that $p_jn\to q_j$ for every $j$ and some $q_j\in[0,\infty)$. (Thus, Theorem 8.2 is the case when all $q_j=0$.) If we further assume, without loss of generality, as in Theorem 8.2, that $\frac{1}{2}\sum_j(p_jn)^2\to\lambda$ for some $\lambda<\infty$, and let $\lambda':=\lambda-\frac{1}{2}\sum_j q_j^2$, it can be shown by arguments similar to those above that

$$n-Z_n\xrightarrow{\mathrm{d}}W+\sum_1^\infty\tilde{V}_j,$$

where $W\sim\mathrm{Po}(\lambda')$, $\tilde{V}_j:=V_j-\mathbf{1}_{\{V_j\ge 1\}}$ with $V_j\sim\mathrm{Po}(q_j)$, and all terms are independent. Note that the limit depends on the sequence $\{q_j\}$; thus, in general, different subsequences may converge to different limits, even if the limit $\lambda$ exists.

Similarly, in case (ii), we may rearrange $(p_j)$ into two (finite or infinite) sequences $(p'_j)$ and $(p''_j)$ with $1/n\ge p'_1\ge p'_2\ge\cdots$ and $1/n<p''_1\le p''_2\le\cdots$, and by selecting a subsequence we may assume that $p'_jn\to q'_j$ and $p''_jn\to q''_j$ for some $q'_j\in[0,1]$ and $q''_j\in[1,\infty]$. (If the sequences are finite, extend them by 0's or $\infty$'s.) It can be shown that if $\lambda_1$, $\lambda_2$ and $m$ are as in Theorem 8.3, then

$$Z_n-m\xrightarrow{\mathrm{d}}W'+\sum_j V'_j-W''-\sum_j V''_j,$$

where $W'\sim\mathrm{Po}(\lambda_1-\sum_j q'_j)$, $W''\sim\mathrm{Po}(\lambda_2-\sum_j e^{-q''_j})$, $V'_j\sim\mathrm{Be}(1-e^{-q'_j})$, $V''_j\sim\mathrm{Be}(e^{-q''_j})$, and all terms are independent. We leave the details to the reader.

**9. Fixed distribution.** We briefly discuss a few characteristic properties for the case when the distribution $F$ is kept fixed while $n \to \infty$. We may as in Remark 3.1 assume that the sequence $(Z_n)$ is obtained by throwing balls one after another; thus $Z_1 \le Z_2 \le \cdots$.

Let $M := \#\{j : p_j > 0\}$, the number of distinct values that $X_i$ can take with positive probability. If $M$ is finite, then a.s. all these values are sooner or later assumed by some $X_i$, and thus $Z_n = M$ for large enough $n$. In other words, then $Z_n = M$ whp, and $Z_n \xrightarrow{\text{p}} M$ as $n \to \infty$.

We will therefore in this section assume that $M = \infty$. It is then easily seen that $Z_n \to \infty$ a.s. as $n \to \infty$; similarly $Z(\lambda) \to \infty$ a.s. as $\lambda \to \infty$. Consequently, $\mathbb{E}(Z_n) \to \infty$ as $n \to \infty$ and $\mu(\lambda) = \mathbb{E}(Z(\lambda)) \to \infty$ as $\lambda \to \infty$.

On the other hand, by (2.2) and the dominated convergence theorem

$$\mu(x)/x = \sum_j (1 - e^{-p_j x})/x \to 0 \qquad \text{as } x \to \infty,$$

since $0 \le (1 - e^{-p_j x})/x \le p_j$ and $\sum_j p_j < \infty$; see also Karlin [20] for an alternative proof. In other words, $\mu(x) = o(x)$ as $x \to \infty$, and thus, by Theorem 2.3, $\mathbb{E}(Z_n) = o(n)$ as $n \to \infty$.

Similarly, the $O(1)$ terms in Theorem 2.3 can be improved to $o(1)$; these remainder terms are given in our proof in Section 4 as sums, where each term tends to 0 and domination is provided by the estimates given in our proof.

Finally, $\sum_{p_j x > 1} p_j \to \sum_{p_j > 0} p_j = 1$ as $x \to \infty$; thus we always have $\sum_{p_j x > 1} p_j \ge 1/2$ for large $x$. Hence $\sigma^2(n) \asymp v(n)$ and, for limit results, we only have to consider the case in Section 5.

**Acknowledgment.** We thank an anonymous referee for careful reading and useful comments, and in particular for improving our previous versions of Proposition 5.6 and Lemma 4.2.

Institute of Statistical Science
Academia Sinica
Taipei 115
Taiwan

Department of Mathematics
Uppsala University
PO Box 480, SE-751 06
Uppsala
Sweden
E-mail: svante.janson@math.uu.se